\documentclass[11pt,a4paper]{article}
\usepackage{geometry}
\usepackage{amsfonts}
\usepackage{amssymb , amsmath, amsthm, amsopn, amstext, amscd}
\usepackage{latexsym, enumerate}
\usepackage{graphics}
\usepackage{graphicx}
\usepackage[french,ruled]{algorithm2e}
\usepackage{a4wide}

\newtheorem{prp}{Proposition}

\newtheorem{rmq}{Remark}
\newtheorem{ass}{Assumption}

\DeclareMathOperator*{\argmin}{argmin}

\newcommand{\eps}{\varepsilon}

\title{Sparse Conformal Predictors}
\author{\noindent{Mohamed Hebiri\footnote{hebiri@math.jussieu.fr}} \\ \\
\small{Laboratoire de Probabilit\'{e}s et Mod\`{e}les Al\'{e}atoires, CNRS-UMR 7599,}\\ 
\small{Universit\'{e} Paris 7 - Diderot, UFR de Math\'{e}matiques,}\\ 
\small{175 rue de Chevaleret F-75013 Paris, France.}}

\date{}

\begin{document}

\maketitle

\begin{abstract}
Conformal predictors, introduced by Vovk et al. [16], serve to build prediction intervals by exploiting a notion of conformity of the new data point with previously observed data. In the present paper, we propose a novel method for constructing prediction intervals for the response variable in multivariate linear models. The main emphasis is on sparse linear models, where only few of the covariates have significant influence on the response variable even if their number is very large. Our approach is based on combining the principle of conformal prediction with the $\ell_1$ penalized least squares estimator (LASSO). The resulting confidence set depends on a parameter $\varepsilon>0$ and has a coverage probability larger than or equal to $1-\varepsilon$. The numerical experiments reported in the paper show that the length of the confidence set is small. Furthermore, as a by-product of the proposed approach, we provide a data-driven procedure for choosing the LASSO penalty. The selection power of the method is illustrated on simulated data.
\end{abstract}

\section{Introduction}\label{Sec:Intro}

Consider observations $(x_{i},y_{i})\in\mathbb{R}^{p}\times \mathbb{R}$ for $i \ge 1$ from a linear regression model $y_i = x_{i}' \beta + \xi_i$, where $\beta \in \mathbb{R}^{p}$ is the unknown parameter and the $\xi_i$'s are the noise variables. Suppose we have already collected the dataset $\mathcal{E}_{n}=\left((x_{1},y_{1}),\ldots,(x_{n-1},y_{n-1}),x_{new}\right)$ where $x_{new}\in \mathbb{R}^{p}$ denotes a new observation. Our goal is to predict the label $y_{new}$ corresponding to $x_{new}$ based on $\mathcal{E}_{n}$ and then exploiting the information in $x_{new}$. This setup is known as the transduction problem~\cite{VapnikLivre}. Our estimation strategy is based on local arguments in order to produce a better estimation for $y_{new}$~\cite{localEstimation}. More precisely, we will follow the approach of {\em conformal prediction} presented by Vovk~et~al.~\cite{Vovk_livre} which relies on two key ideas: one is to provide a confidence prediction (namely, a confidence set containing $y_{new}$ with high probability) and the other is to account for the similarity of the new data $x_{new}$ compared to the previously observed $x_i$'s. The notion of conformal predictor was first described by Vovk et al.~\cite{Vovk99IntroCP}. Moreover, in~\cite{Vovk_livre}, the authors illustrate this approach on the example of ridge regression. Along the paper, this predictor will be referred to as Conformal Ridge Predictor\footnote{The Conformal Ridge Predictor was called the Ridge Regression Confidence Machine in Vovk et al. \cite{Vovk_livre}.} (CoRP). In the present contribution, we propose to adapt conformal predictors to the sparse linear regression model, that is a model where the regression vector $\beta\in\mathbb{R}^{p}$ contains only a few of nonzero components. We introduce a novel conformal predictor called the {\it Conformal Lasso Predictor} (CoLP) which takes into account the sparsity of the model. Its construction is based on the LASSO estimator \cite{Tibshirani-LASSO}. The LASSO estimator for linear regression corresponds to an $\ell_1$-penalized least square estimator and it has been extensively studied over the last few years (\cite{Knight&Fu-Lasso-distrib, Meins-Buhl-Lasso-Graph, TsybSparsLasso, Zhao-Yu-Consist-Lasso}, among others) and
several modifications have been proposed (\cite{Zou-Adaptive-Lasso, Yuan-Lin-GroupSelection, Zou-E-Net, Rosset-Fused, MoiSLasso} among others). One attractive aspect of the LASSO is that it aims both to provide accurate estimating while enjoying variable selection when the model is sparse. In the approach considered in the present paper, the resulting Conformal Lasso Predictor has a large coverage probability and are small in term of its length in the same time. When we deal with regularized methods like the Ridge or the LASSO estimators, the choice of the penalty is an important task. Contrary to the Conformal Ridge Predictor for which no rule was established to pick the Ridge-penalty~\cite{Vovk_livre}, the construction of the Conformal Lasso Predictor provides a data-driven way for choosing the LASSO-penalty. Moreover, it turn out that this choice is adapted to variable selection as supported by the numerical experiments.

The paper is organized as follows. We concisely introduce conformal prediction and the LASSO procedure in Section~\ref{sec:ConfidencePredictor} and Section~\ref{LassoOriginal} respectively. In Section~\ref{section-C-Lasso}, we give the explicit form of the Conformal Lasso Predictor. An algorithm producing the CoLP is presented in Section~\ref{sec:Algo}. Then in Section~\ref{sec:generalization} we discuss a generalization of the Conformal Lasso Predictor to other selection-type procedures; we call these generalized procedures {\it Sparse Conformal Predictors}. Finally, in Section~\ref{sec:simul}, we illustrate the performance of Sparse Conformal Predictors through some numerical experiments.

\section{Conformal prediction}\label{sec:ConfidencePredictor}
Let us briefly describe the approach based on conformal prediction developed in the book by Vovk et al. \cite{Vovk_livre} where they develop the idea of {\it conformal} prediction. In order to predict the  label $y_{new}$ of a new observation $x_n = x_{new}$, the similarity of pairs of the form $(x_{new},y)$, where $y\in\mathbb{R}$,  to the former observations $(x_{i},y_{i})$ for $i=1,\ldots,n-1$ is exploited. This is the purpose of introducing a {\it nonconformity score} $\alpha(y) = (\alpha_{1}(y),\ldots,\alpha_{n}(y))'$ which is based on $\mathcal{E}_{n}$. Each component $\alpha_{i}$ describes the efficiency of explaining the observation $(x_i,y_i)$ by a procedure based on the augmented sample $\left\{(x_{1},y_{1}),\ldots,(x_{n-1},y_{n-1}),(x_{new},y)\right\}$. In order to obtain a relative information between  different nonconformity scores $\alpha_{i}$, we shall use the notion of {\it $p$-value}, as introduced in~\cite{Vovk_livre}, defined as:
\begin{equation}\label{p-value1}
p(y) = \frac{1}{n}|\left\{i\in\{1,\ldots,n\}:\, \alpha_{i}(y)\geq
\alpha_{n}(y)\right\}|~,
\end{equation}
where for any set $\mathcal{A}$, we denote its cardinality by $|\mathcal{A}|$. The above quantity lies between $1/n$ and $1$. Moreover, we note that the smaller this $p$-value is, the less likely the tested pair $(x_{new},y)$ is (in other words, $y$ is an outlier when associated to $x_{new}$). An explicit form of the nonconformity score and the $p$-value will be given in Section~\ref{section-C-Lasso} when we will adapt it to the CoLP.
\begin{rmq}
The notion of $p$-value introduced in the present paper differs from the classical one. To make the connection with hypothesis testing in mathematical statistics~\cite{pValueCasellaBerger}, consider the following hypotheses:
\[ \left \lbrace 
\begin{array}{cl}
H_0: & \text{the pair $(x_{new},y)$ is conformal}, \\
H_1: & \text{the pair $(x_{new},y)$ is not conformal.}
\end{array}
\right. \]
Assume the observation $Y=y$ is given. The function $p(y)$ permits to construct a statistical test procedure with critical region $\mathcal{R}_{\varepsilon} = \{y: \, p(y) \leq \varepsilon \}$ and $H_0$ is rejected if $y\in \mathcal{R}_{\varepsilon}  $.
\end{rmq}
A nice feature of this nonconformity score is that it can be related to the confidence
of the prediction for $y_{new}$. We now recall the concept of conformal predictor introduced in
\cite{Vovk_livre}. Set $\varepsilon\in(0,1)$. Given the new observation $x_{new}$, we search for a
subset $\Gamma^{\varepsilon} = \Gamma^{\varepsilon}(\mathcal{E}_{n})$ of
$\mathbb{R}$, in which the expected value of $y_{new}$ lies with
a probability of $1-\varepsilon$. The conformal predictor $\Gamma^{\varepsilon}$ is defined as the set of labels $y\in \mathbb{R}$ such that $p(y)>\varepsilon$. {\it In other words, $\Gamma^{\varepsilon}$ consists of labels $y$ which make the pair $(x_{new},y)$ more conformal than a proportion $\varepsilon$ of the previous pairs $(x_{i},y_i)$ for $i=1,\ldots,n-1$}. Note moreover that the smaller $\varepsilon$, the more confident the predictor. That is to say, for any $\varepsilon_{1}, \,\varepsilon_{2} >0$:
\begin{equation*}
\Gamma^{\varepsilon_{1}}\subset \Gamma^{\varepsilon_{2}}  \quad
\text{whenever}\,\, \varepsilon_{1} \geq \varepsilon_{2}~.
\end{equation*}

In the present analysis, apart from prediction, we develop an approach for selecting relevant variables. For this reason, we consider three criteria measuring the quality of our procedure: {\em validity}, {\em accuracy}, and {\em selection}.
The first two were introduced in \cite{Vovk_ConfLinReg}. The fact that we consider the issue of sparsity leads us to include the selection power of the predictor.
\begin{description}
\item[Validity.]
This criterion accounts for the power of conformal prediction. The simplest approach is to count the number of times where $y_n$ does not belong to the set $\Gamma^{\varepsilon}$. We take the notation:
\begin{equation*}
\text{err}_{n}^{\varepsilon} = \left\{
\begin{array}{c}
1 \quad \text{if}\, y_{n} \notin  \Gamma^{\varepsilon} (\mathcal{E}_n)  \\
0 \quad \text{otherwise.}
\end{array}
\right.
\end{equation*}
Note that in an on-line perspective, one focuses on the cumulative error $\text{ERR}_{n}^{\varepsilon} = \sum_{i=1}^{n}\text{err}_{i}^{\varepsilon}$. Asymptotic validity properties of this cumulative error have been studied in~\cite{Vovk02AsymptVal} and~\cite[chapters 2 and 8]{Vovk_livre}. In the present work, we will be interested in evaluating the error $\text{err}_{n}^{\varepsilon}$ for a fixed $n$, rather than the cumulative one.\\
\item[Accuracy.] The length of the confidence predictor provides a natural measure of the accuracy. We will see that such a measure is adapted to the variable selection purpose. Note that other choices are possible. We shall discuss this point in Section~\ref{sec:Algo}.\\
\item[Selection.] Finally, in the case of sparse linear regression, it is important to include a measure of
the capacity of the estimator to select relevant variables, namely those for which the regression parameter $\beta$ has nonzero components.
\end{description}

\section{The LASSO Procedure}\label{LassoOriginal}The LASSO estimator \cite{Tibshirani-LASSO} has originally been introduced in
the linear regression model:
\begin{equation} \label{eq_depart}
y_{i}= x_{i}' \beta^{*} + \xi_{i}, \quad \quad i=1,\ldots,n-1
\end{equation}where the design $x_{i}=(x_{i,1},\ldots,x_{i,p})' \in \mathbb{R}^p$ is
deterministic, $\beta^*=(\beta^*_1,\ldots,\beta^*_p)' \in \mathbb{R}^p$ is the
unknown regression vector and the $\xi_i$'s are independent and
identically distributed (i.i.d.) centered Gaussian random variables with known
variance $\sigma^{2}$. Then the goal is to use the observations to
provide an approximation of the label $y_{new}$ of a new observation $x_{new}$
through the estimation of the regression vector $\beta^{*}$. The LASSO
estimator is defined as follows:
\begin{equation}\label{lasso-criter}
\hat{\beta}_{\lambda} =
\argmin_{\beta\in\mathbb{R}^{p}}\sum_{i=1}^{n-1}\left(y_{i}-
x_{i}'\beta\right)^{2}+ \lambda\sum_{j=1}^{p}|\beta_{j}|,
\end{equation}where $\lambda\ge 0$ is a tuning parameter. Based on $\hat{\beta}_{\lambda}$, an estimation of the response $y_{new}$ of the new observation $x_n = x_{new}$ is produced by $\hat{\mu}_{\lambda}=x_{new}' \hat{\beta}_{\lambda}$. For a large enough $\lambda$, the LASSO estimator is sparse. That is many components of $\hat{\beta}_{\lambda}$ equal zero. Therefore we can naturally define a sparsity (or active) set as $\mathcal{A}_{\lambda} = \{j\in\{1,\ldots, p\}:\,\hat{\beta}_{\lambda}\neq 0\}$. A LASSO modification of the LARS algorithm~\cite{Efron-LARS} can iteratively provide approximations of the LASSO estimator for a few values of the tuning parameters $\lambda = \lambda_{0},\ldots,\lambda_{K}$ such that $\infty=\lambda_{0}>\ldots>\lambda_{K}=0$ (the indices refer to the algorithm steps and $K$ denotes the last step). These points are the so-called {\it transition points}.\\
From now on, let us write $\hat{\beta}_{k}$ and $\mathcal{A}_k$ for the LASSO estimator $\hat{\beta}_{\lambda}$ and the sparsity set $\mathcal{A}_{\lambda}$ evaluated at the transition point $\lambda = \lambda_{k}$. Obviously, the estimator $\hat{\beta}_{k}$ is an $|\mathcal{A}_{k}|$-dimensional vector where $|\mathcal{A}_{k}|$ is the cardinality of the set $\mathcal{A}_{k}$. Furthermore, we denote by $s_{k}$ the $|\mathcal{A}_{k}|$-dimensional sign vector whose components are the signs of the components of the LASSO estimator evaluated at the transition point $\lambda_k$ (i.e., $(s_{k})_j=1$ if $(\hat{\beta}_{k})_j > 0$, $(s_{k})_j=-1$ if $(\hat{\beta}_{k})_j < 0$ where $j\in\mathcal{A}_{k}$). Finally, let us denote by $\mathbf{x}_{k}$, the $(n-1)\times|\mathcal{A}_{k}|$ matrix whose columns are the variables $X_{j}$, with indices $j\in\mathcal{A}_{k}$. For each $\lambda_k$, we assume that the matrix $(\mathbf{x}_{k}'
\mathbf{x}_{k})^{-1}$ is invertible. Here are some characteristics of the LARS algorithm and we refer to~[2] for more details:\\
\begin{enumerate}
\item[i)] At each iteration of the algorithm (i.e., at each transition point), only one variable $X_{j}=(x_{1,j},\ldots,x_{n-1,j})',\,j=1,\ldots,p$ is added (or deleted) to the construction of the estimator according to its correlation with the current residual. The algorithm begins with only one variable and ends up with the ordinary least square (OLS) estimator\footnote{When $p>n$, the LARS cannot select all $p$ variables. It is limited by the sample size $n$. In such a case, the last iteration does not correspond to the OLS.}.\\
\item[ii)] For each $\lambda\in\left(\lambda_{k+1},\lambda_{k}\right]$, the LASSO estimator can be expressed in the following form:
\begin{equation}\label{eq:PetitLasso}
\hat{\beta}_{\lambda}(\mathbf{y},\mathbf{x}_{k},s_{k})=(\mathbf{x}_{k}'
\mathbf{x}_{k})^{-1} (\mathbf{x}_{k}'\mathbf{y} - \frac{\lambda}{2} \, s_{k} ),
\end{equation}
where $\mathbf{y} =\left(y_{1},\ldots,y_{n-1}\right)'$. Note that \eqref{eq:PetitLasso} is obtained by minimizing \eqref{lasso-criter} over the set $\mathcal{A}_{k}$. Let us also mention that the set $\mathcal{A}_{k}$ and the sign vector $s_{k}$ remain unchanged when $\lambda$ varies in the interval $\left(\lambda_{k+1},\lambda_{k}\right]$.\\
\item[iii)] As highlighted by \eqref{eq:PetitLasso}, the LASSO estimator is piecewise linear in $\lambda$ and linear in $\mathbf{y}$ for every fixed $\lambda$~\cite{Rosset-PeacewiseLin}. Using the LASSO modification of the LARS algorithm, this property helps us to provide the regularization path of the LASSO estimator, which is defined as $\{\hat{\beta}_{\lambda}:\,\, \lambda\in \left[0,\infty\right)\}$ (each point of the regularization path corresponds to the evaluation of the regression vector estimator for a given value of $\lambda$). Indeed, the slope of the LASSO regularization path changes at a finite number of points which coincide with the transition points $\lambda_{1},\ldots,\lambda_{K}$.\\
\item[iv)] Piecewise linearity is an important property of the LASSO modification of the LARS algorithm. Indeed, let $\lambda\in\left(\lambda_{k+1},\lambda_{k}\right]$ where $\lambda_{k+1}$ and $\lambda_{k}$ are two transition points. In this interval, the LASSO estimator  $\hat{\beta}_{\lambda}$ uses the same variables (variables with indices in $\mathcal{A}_{k}$). By using \eqref{eq:PetitLasso}, it is easy to see  \cite{Zou-df-Lasso} that the linearity of the LASSO estimator implies that, for any $\lambda\in\left(\lambda_{k+1},\lambda_{k}\right]$:
\begin{equation*}
\sum_{i=1}^{n-1}\left(y_{i}-x_{i}'\hat{\beta}_{\lambda}\right)^{2}>
\sum_{i=1}^{n-1}\left(y_{i}-x_{i}'\hat{\beta}_{\lambda_{k+1}}\right)^{2}.
\end{equation*}
This last observation indicates that the transition points are the most interesting points in the regularization path.
\end{enumerate}
All these nice properties encourage the use of the LASSO as a selection procedure. In the sequel, we will consider the LASSO modification of the LARS algorithm which provides an approximate solution to the LASSO.
\begin{rmq}
Through the paper, one should keep in mind the analogy between each iteration $k$ of the modification of the LARS algorithm and its corresponding tuning parameter value $\lambda_k$. Decrease of tuning parameter $\lambda$ is reflected through the increase of the number of iterations of the modification of the LARS algorithm.
\end{rmq}
\section{Sparse predictor with conformal Lasso}\label{section-C-Lasso}

For the reasons exposed above, we focus on the transition points $\lambda_{1},\dots,\lambda_{K}$ and construct conformal predictors for each of these $\lambda_{k}$. We then propose to select the best conformal predictor among them according to its performance in terms of accuracy (cf. Section~\ref{sec:ConfidencePredictor}).

Now let us detail the construction of the CoLP for each $\lambda_k$. To this end, denote by $X_j=(x_{1,j},\ldots,x_{n-1,j},x_{new,j})',\,j=1,\ldots,p$ the augmented variable $j$. Define the augmented matrix $\widetilde{\mathbf{x}}=\left(x_{1},\ldots,x_{n-1},x_{new}\right)'= \left(X_{1},\ldots,X_{p}\right)$ and the augmented response vector $\widetilde{\mathbf{y}}=\left(y_{1},\ldots,y_{n-1},y\right)'$ where $y$ is a candidate value for $y_{new}$. Using the notation introduced in Section~\ref{LassoOriginal}, for the fixed $\lambda_k$, we also define the LASSO estimator $\hat{\beta}_{k}(\widetilde{\mathbf{y}},\widetilde{\mathbf{x}}_k,s_{k})$ from expression \eqref{eq:PetitLasso} with the augmented data. From now on, we denote this estimator by $\hat{\beta}_{k}$. Define $\hat{\mu}_{k}:=\widetilde{\mathbf{x}}_{k} \hat{\beta}_{k} $. Moreover, the matrix $\mathbf{H}_{k}$ will be the $n\times n$ projection matrix onto the subspace generated by $\widetilde{\mathbf{x}}_{k}$ and $\mathbf{I}$ identity matrix of the same size. For each $\lambda_k$, we define a corresponding nonconformity score $\alpha^{k}=\left(\alpha_{1}^{k},\ldots,\alpha_{n}^{k}\right)'$ by:
\begin{eqnarray*}
\alpha^{k}(y) & := & |\widetilde{\mathbf{y}}-\hat{\mu}_{k}| =
|\left(\mathbf{I}-\mathbf{H}_{k}\right)\widetilde{\mathbf{y}} +
\frac{\lambda_{k}}{2} \widetilde{\mathbf{x}}_{k} \left(\widetilde{\mathbf{x}}_{k}'
\widetilde{\mathbf{x}}_{k}\right)^{-1} s_{k}|
\\ & = &|A_{k} + B_{k}\,y|,
\end{eqnarray*}
where $|\cdot|$ is meant here componentwise and
\begin{equation}\label{eq:AkBkCk}
\left\{
\begin{array}{l}
A_{k}= (a_{1}^{k},\ldots,a_{n}^{k})' :=\left(\mathbf{I}-\mathbf{H}_{k}\right) (y_{1},\ldots,y_{n-1},0)' + 
\frac{\lambda_k}{2} \widetilde{\mathbf{x}}_{k}
\left(\widetilde{\mathbf{x}}_{k}' \widetilde{\mathbf{x}}_{k}\right)^{-1} s_{k}, \\
B_{k}= (b_{1}^{k},\ldots,b_{n}^{k})' :=\left(\mathbf{I}-\mathbf{H}_{k}\right) (0,\ldots,0,1)',\\
\end{array}
\right.
\end{equation}
Note that each component $\alpha_{i}^{k}(y)$ is piecewise linear with respect to $y$. Then the corresponding $p$-value $p_{k}(y)$ as defined by \eqref{p-value1} clearly can change only at points $y$ where the sign of $\alpha_{i}^{k}(y)-\alpha_{n}^{k}(y)$ changes. Hence, we do not have to evaluate all the possible values of $y$. We only focus on points $y$ for which the $i$-th nonconformity measure $\alpha_{i}^{k}(y)$ equals $\alpha_{n}^{k}(y)$. For this purpose, we define, for each observation $i\in\{1,\ldots,n\}$
\begin{equation}\label{equat-S}
S_{i}^{k}=\left\{y:\,\alpha_{i}^{k}(y)\geq \alpha_{n}^{k}(y)\right\},
\end{equation}
which corresponds to the range of values $y$ such that the new pair $(x_{new},y)$ has a better conformity score than the $i$-th pair $(x_{i},y_{i})$. Moreover, let $l_i^k$ and $u_i^k$ denote two real defined respectively as
\begin{equation}\label{UpandDown}
l_i^k = \min\{-\frac{a_{i}^{k}-a_{n}^{k}}{b_{i}^{k}-b_{n}^{k}};-\frac{a_{i}^{k}+a_{n}^{k}}{b_{i}^{k}+b_{n}^{k}}\}, \quad\quad\text{and }\quad u_i^k = \max\{ -\frac{a_{i}^{k}-a_{n}^{k}}{b_{i}^{k}-b_{n}^{k}};-\frac{a_{i}^{k}+a_{n}^{k}}{b_{i}^{k}+b_{n}^{k}}\},
\end{equation}
where $a_i^k$ and $b_i^k$ are given by~\eqref{eq:AkBkCk}.

\begin{prp}\label{pr:concide}
Let us fix a $k\in\{1,\ldots,K\}$ and an $i\in\{1,\ldots,n-1\}$. Assume that both $b_i^k$ and $b_n^k$ are non-negative. Then\\
\begin{itemize}
\item[i)] if $b_{i}^{k}\neq b_{n}^{k}$, we have either $S_i^k = [l_i^k ; u_i^k]$ or $S_i^k = (-\infty ; l_i^k] \cup [u_i^k ; -\infty )$, with $l_i^k$ and $u_i^k$ given by~\eqref{UpandDown}.\\
\item[ii)] if $b_{i}^{k}=b_{n}^{k}\neq 0$, then $l_i^k = u_i^k = -\frac{a_{i}^{k}+a_{n}^{k}}{2 b_{n}^{k}}$ and we have either $S_i^k = (-\infty ; l_i^k]$ or $S_i^k = [l_i^k ; -\infty )$. Moreover if $a_{i}^{k}=a_{n}^{k}$, we have $S_i^k = \mathbb{R}$.\\
\item[iii)] if $b_{i}^{k}=b_{n}^{k}=0$, we have either $S_i^k = \mathbb{R}$ or $S_i^k= \emptyset$.
\end{itemize}
\end{prp}
The assumption that all the $b_i^k $ are non-negative does not make loose any generality as one can multiply $a_{i}^{k}$, $b_{i}^{k}$ and $c_{i}^{k}$ by $-1$ if $b_{i}^{k}<0$. With this definition of $S_i^k$, we may rewrite the definition of the conformal predictor as follows
\begin{equation}\label{eq:ConfPredick}
\Gamma_k^{\varepsilon}=\{y:\sum_{i=1}^n \mathbb{I}(\alpha_i^{k}(y)\ge \alpha_n^k(y))\ge n\varepsilon\}=\{y:\sum_{i=1}^n \mathbb{I}({S^k_i})(y) \ge n\varepsilon\},
\end{equation}
where $\mathbb{I}(\cdot)$ stands for the indicator function. This approach leads to a whole collection of confidence intervals 
$\Gamma^{\eps}_{1},\ldots,\Gamma^{\eps}_{K}$. We propose below a strategy for choosing one
one particular $\Gamma^\eps_k$, the performance of which will be studied through numerical simulations. 

It is worth mentioning that in view of~\cite[Theorem~1]{Vovk02IndepError} (see also~\cite[Proposition~2.3 page~26]{Vovk_livre}), 
each of predictor $\Gamma^\eps_k$ would have a coverage probability at least equal to $1-\varepsilon$, if the corresponding value $\lambda_k$
of the tuning parameter were deterministic. In fact, the following result holds.
\begin{prp}\label{prp:Val}
Fix the significance level $\varepsilon\in (0,1)$ and the tuning parameter $\lambda>0$. Let $\hat\beta_{\lambda,n}(y)$ be the Lasso estimate for the augmented dataset $(\tilde{\mathbf y},\tilde{\mathbf x})$ and let us define $\alpha^\lambda(y)=
|\tilde{\mathbf y}-\tilde{\mathbf x}\hat\beta_{\lambda,n}(y)|$. Then, the conformal predictor 
$$
\Gamma_{\lambda}^{\varepsilon}=\Big\{y:\sum_{i=1}^n \mathbb{I}(\alpha_i^{\lambda}(y)\ge \alpha_n^{\lambda}(y))\ge n\varepsilon\Big\},
$$ 
satisfies
\begin{equation*}
\mathbb{P}(y_{new}\in \Gamma_k^{\varepsilon}) \geq 1-\varepsilon,
\end{equation*}
for any $n\in \mathbb{N}$.
\end{prp}

Actually, in the proof of Proposition~\ref{prp:Val} detailed in~\cite{Vovk02IndepError}, one needs the exchangeability of the pairs $(x_1,y_1),\ldots,(x_{n-1},y_{n-1}),(x_n,y)$ in the definition of the predictor. 
This property is not fulfilled when the tuning parameter $\lambda$ is chosen in the set $\{\lambda_1,\ldots,\lambda_K\}$ of Lasso's transition points, since the elements of this set depend only on the first $n-1$ 
observations and not on 
$(x_n,y)$. We believe that under some additional assumptions a result similar to Proposition~\ref{prp:Val} can be obtained for the predictor $\Gamma^\eps_k$ as well, for each $k=1,\ldots,K$. This is the topic of an 
ongoing work. In the present paper, we content ourselves by proposing a data-driven choice of the conformal predictor from the collection of predictors $\{\Gamma_k^\eps;1\le k\le K\}$ and by exploring its empirical 
properties. 
\begin{rmq}
Of course, one can also apply the well-known sample splitting technique for choosing the values $\lambda_1,\ldots,\lambda_K$ based on a first sample, and then use the methodology described below
for selecting the data-driven predictor based on a second sample which is assumed to be independent of the first sample. However, this technique is not attractive from the practical standpoint, that is why
we do not develop this approach.
\end{rmq}
As discussed above, we believe that all the predictors $\Gamma_k^\epsilon$ share nearly the $1-\varepsilon$ validity property, which is supported by our empirical study. We suggest to select among them the one which has the smallest Lebesgue measure. We denote this confidence set by $\Gamma^\epsilon_{opt}$, that is
\begin{equation}\label{eq:CoLPoptimal}
 \Gamma^{\varepsilon}_{opt}=\Gamma^{\eps}_{\nu},\qquad \nu=\argmin_{k} |\Gamma^{\varepsilon}_{k}|.
\end{equation}
In general, since $\nu$ is a random variable, the $1-\varepsilon$ validity of all $\Gamma^{\varepsilon}_{k}$ would not imply the $1-\varepsilon$ validity of $\Gamma^{k}_{opt}$, but only $1-K\varepsilon$ validity. However, $1-K\eps$ is a worst case majorant obtained by a simple application of the union bound, whereas numerical examples we considered (some of them are reported below) suggest that the validity is much better than $1-K\varepsilon$ and could even be equal to $1-\varepsilon$ when~$p\leq n$. 
\section{Implementation}\label{sec:Algo}
We provide here a three-step algorithm which enables us to easily construct the CoLP. We start in \textbf{Step~1} by applying the LASSO modification of the LARS algorithm to the dataset $\left((x_{1},y_{1}),\ldots,(x_{n-1},y_{n-1})\right)$. This step provides all transition points $\lambda_{1},\ldots,\lambda_{K}$, the corresponding design matrices $\mathbf{x}_{k}$ and sign vectors $s_{k}$ for $k=1,\ldots,K$. Then, in \textbf{Step~2}, we construct the conformal predictor $\Gamma_{k}^{\varepsilon}$ associated to each $\lambda_{k}$. Thanks to Proposition~\ref{pr:concide}, for each $\lambda_{k}$, we can construct the sets $S_{i}^{k}$ for $i=1\ldots,n$ defined by \eqref{equat-S}. We use these sets in order to construct the conformal predictor $\Gamma_{k}^{\varepsilon}$. To do this, we take advantage from the fact that the function $y\mapsto \sum_{i=1}^n \mathbb{I}({S_i^k}(y))$ is piecewise constant. Furthermore, the endpoints of the intervals where this function is constant belong to the set of the all endpoints of intervals forming the sets $S_i^k$. Thus, to determine $\Gamma^\eps_k$, we sort the set $U$ consisting of the all endpoints of the intervals described in Proposition 1 and include an interval having as endpoints two successive elements of $U$ in $\Gamma_k^\eps$ if the center of this interval belongs to at least $[n \eps]  $ sets $S_i^k$.

\begin{figure}[t]
\vskip -0.2in
\includegraphics[height=2.3in] {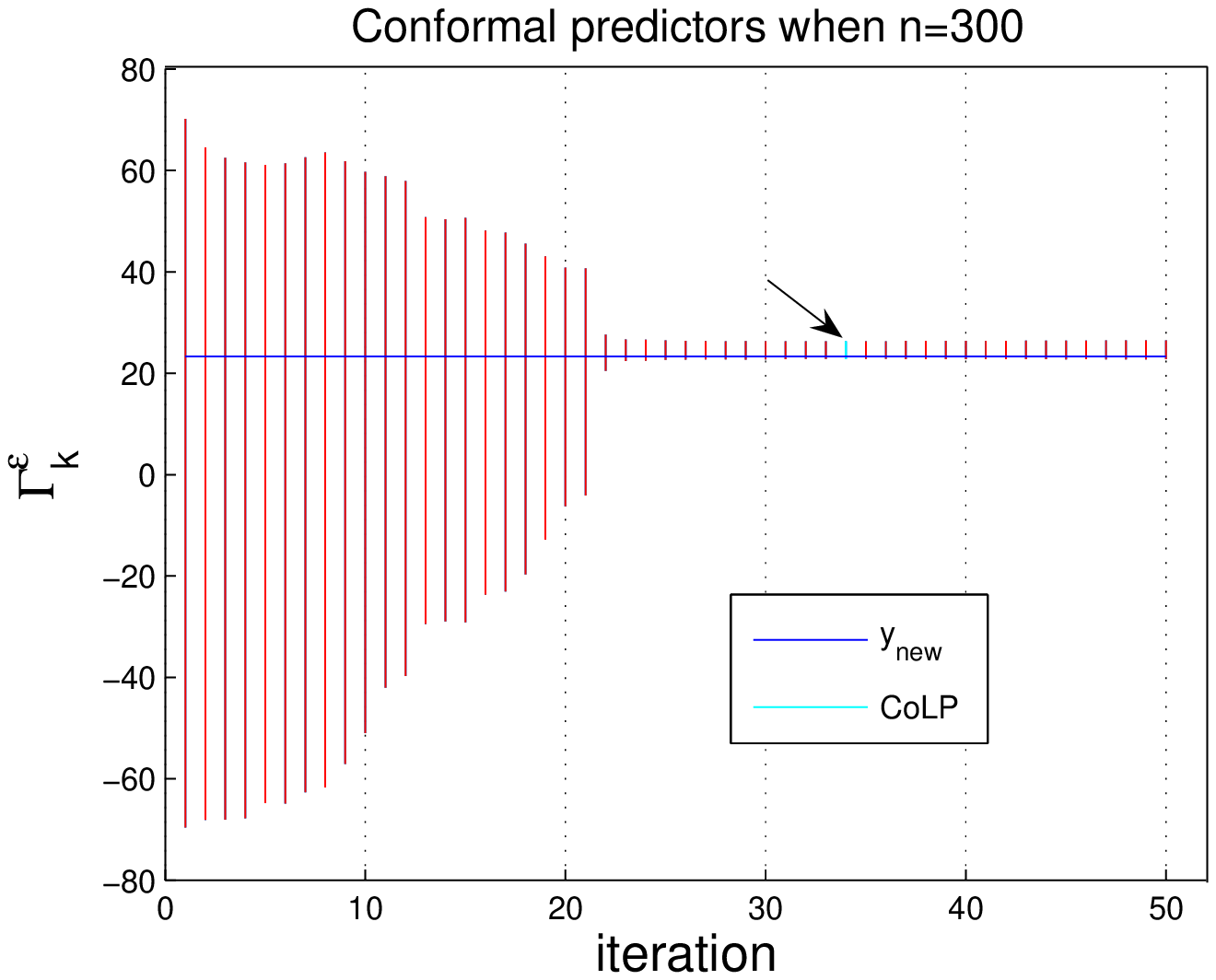}
\includegraphics[height=2.3in] {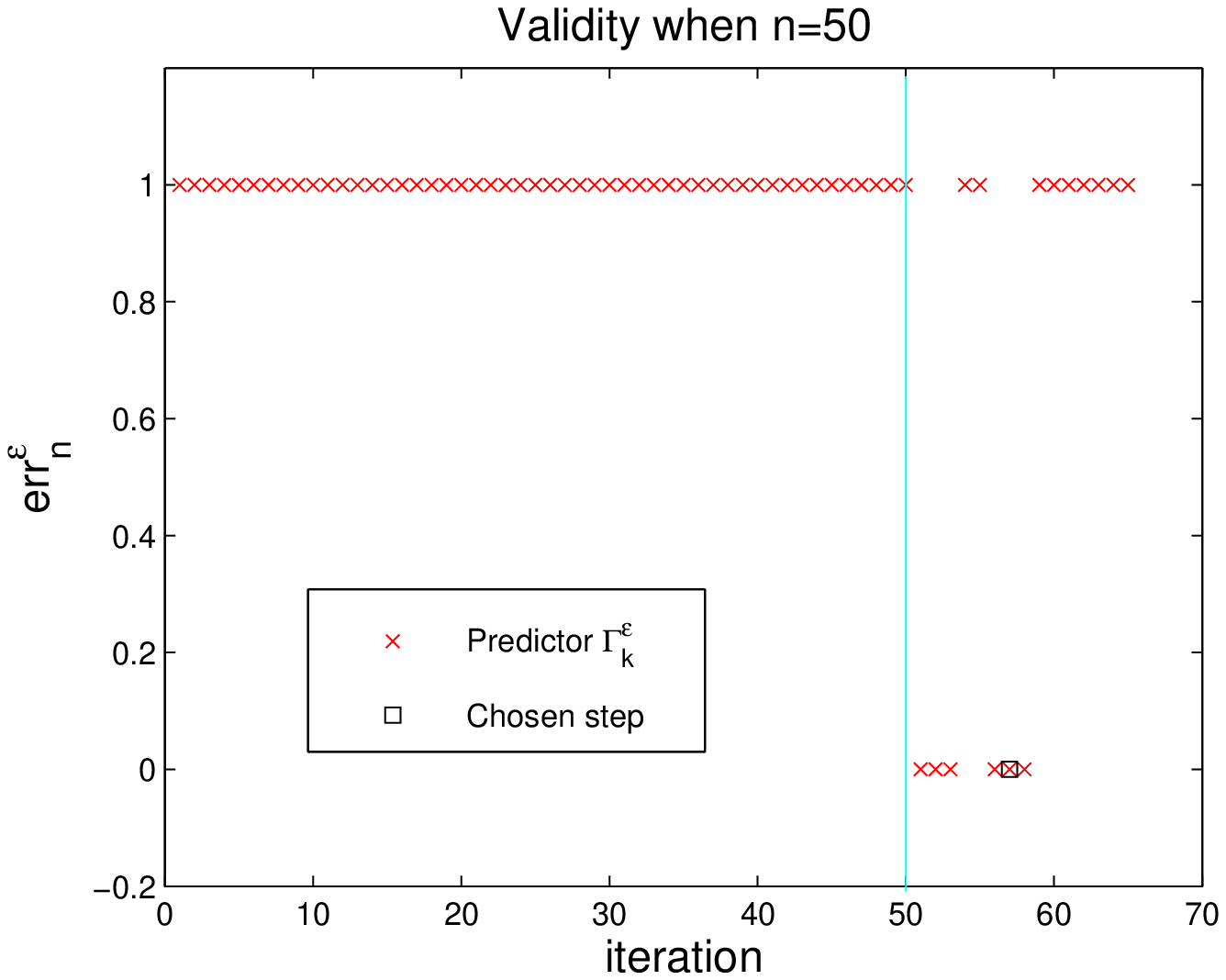}
\vglue-40pt
~
\begin{center}
\begin{minipage}[t]{0.90\textwidth}
\caption{\footnotesize{\it Left}: Conformal predictors $\Gamma_k^{\varepsilon}$ evolution through the iterations of the LASSO modification of the LARS algorithm when $n=300$ (the first iteration corresponds to $\lambda_{max}$ and the last one corresponds to $\lambda_{min}$). The CoLP is drawn in cyan and corresponds to the $34$-th iteration. The horizontal blue line corresponds to the value of $y_{new}$. {\it Right}: Validity analysis ($\text{err}_{n}^{\varepsilon}$) of the conformal predictors $\Gamma_k^{\varepsilon}$ through the iterations of the LASSO modification of the LARS algorithm when $n=50$ (the first iteration corresponds to $\lambda_{max}$ and the last one corresponds to $\lambda_{min}$). The CoLP is marked by a black square and corresponds to the $57$-th iteration. The vertical line represents a separation between a stable and an unstable zone.}
\end{minipage}
\end{center}
\label{fig:SuccesVsSteps}
%\vspace{-0.27in}
\end{figure}

\begin{algorithm}[H]
\dontprintsemicolon {\bfseries Step 1:} {Run the LASSO modification of the LARS
algorithm on the data set
$\left((\mathbf{x}_{1},y_{1}),\ldots,(\mathbf{x}_{n-1},y_{n-1})\right)$}\;
{\bfseries Step 2:} {Construct the Conformal Lasso Predictors for each
$\lambda_{k}\in \{\lambda_{1},\ldots,\lambda_{K}\}$} \Begin{ {\bfseries Step
2a:} {Initialization : Define $A_{k}$ and $B_{k}$ as in~\eqref{eq:AkBkCk}. Set $U^{k}\longleftarrow  \emptyset$}\;
{\bfseries Step 2b:} {Harmonization}\; \For{$i=1$ {\bf to} $n$}{ \If{$b_{i}^{k}
< 0$}{$a_{i}^{k}=-a_{i}^{k}$ and $b_{i}^{k}=-b_{i}^{k}$ } } {\bfseries Step 2c:} {Actualize the set $U^{k}$}\;
\For {$i=1$ {\bfseries to} $n$}{ \If{$b_{i}^{k}\neq b_{n}^{k}$}{Add $l_i^k$ and $u_i^k$~\eqref{UpandDown} to $U^{k}$}\; 
\If{$b_{i}^{k} = b_{n}^{k} \neq 0$ and $a_{i}^{k} \neq a_{n}^{k} $
}{Add $l_i^k = u_i^k$~\eqref{UpandDown}
to $U^{k}$}\; } 
{\bfseries Step 2d:} {Sort $U^{k}$.
Let $m\longleftarrow|U^{k}|$. Then $y_{(0)}\longleftarrow-\infty$ and
$y_{(m+1)}\longleftarrow +\infty$}\; 
{\bfseries Step 2e:} {Evaluate $N_{j}^{k}$ for $j=1,\ldots,m$. Initialize $N_{j}^{k} \longleftarrow 0$. Then actualize}\;{
\For {$i=1$ {\bf to} $n$}{
\For {$j=1$ {\bf to} $m$}{\If{$|a_{i}^{k}+ b_{i}^{k}\,y| \geq |a_{n}^{k} +  b_{n}^{k}\,y|$ \text{for}
$y \in (y_{(j)},y_{(j+1)})$}{ Increment $N_{j}^{k}=N_{j}^{k}+1$}\;
}}} 
{\bfseries Step 2f:} {For a fixed threshold $\varepsilon>0$, output the conformal predictor}
\begin{equation*}
     \Gamma_{k}^{\varepsilon}=\cup_{j:\frac{N_{j}}{n}>\varepsilon}
     [y_{(j)},y_{(j+1)}] 
\end{equation*}
} {\bfseries Step 3:} {Output the Conformal Lasso Predictor $\Gamma_{opt}^{\varepsilon}$ as the smallest
(w.r.t. their Lebesgue measure) confidence set among the constructed conformal predictors}
\caption{Lasso Conformal Predictor
   \label{alg:ConformalLasso}}
\end{algorithm}\\
Finally, in a \textbf{Step~3}, we provide the CoLP, says $\Gamma_{opt}^{\varepsilon}$, which is defined as the smallest confidence set, according to its Lebesgue measure, among the constructed conformal predictors $\Gamma_{k}^{\varepsilon},\,k=1,\ldots,K$. According to Proposition~\ref{prp:Val}, each $\Gamma_{k}^{\varepsilon}$ is valid. Moreover the criterion for choosing the CoLP is adapted to variable selection as conformal predictors constructed here for different values of $\lambda_{k},\,k=1,\ldots,K$ bring into play different variables. This is illustrated in Figure~\ref{fig:SuccesVsSteps} (left) where we constructed the conformal predictors when $n=300$. One can observe that all the conformal predictors are valid since they contain the true value of the label $y_{new}$. Hence our construction is suitable when the sample size is larger than the number of variables (i.e., $n>p$) but may be not appropriated when $p\geq n$. Figure~\ref{fig:SuccesVsSteps} (right) shows an example where almost all the constructed conformal predictors $\Gamma_{k}^{\varepsilon},\,k=1,\ldots,K$, using the above algorithm are valid. Only six are not. One of them is the selected CoLP (iteration $57$ in Figure~\ref{fig:SuccesVsSteps} (right)) which corresponds to the smallest predictor. In such cases ($p\geq n$), a correction can be made and other choices for the accuracy measure are possible. We discuss this criterion in Section~\ref{sec:simul}. Let us add that we only illustrated the validity of the conformal predictors in Figure~\ref{fig:SuccesVsSteps} (right) as the unstable zone (on the right side of the vertical line) makes the representation hard to be  analyzed. More details are given in Section~\ref{sec:simul}.
\begin{rmq}\label{rq:online}
In \textbf{Step~1} of Algorithm~\ref{alg:ConformalLasso}, we use the LARS algorithm for its ability to generate a small number of tuning parameter values of interest. It is an important aspect as it considerably reduces the computational cost. On-line versions could be implemented by plugging in an on-line version of the LASSO solution as in~\cite{OnlineLasso}. The analysis of such on-line versions is the object of work under progress.
\end{rmq}

\section{Extension to others procedures}\label{sec:generalization}
In this section we generalize the construction of the confidence predictor to a family of estimators which includes selection-type procedures as the Elastic-Net~\cite{Zou-E-Net} and the Smooth-Lasso~\cite{MoiSLasso}. As for CoLP (Section~\ref{section-C-Lasso}), we are interested in two properties of estimators: the {\it piecewise linearity w.r.t. the response $y$} (to easily compute the nonconformity scores $\alpha_i ,\, i=1,\ldots,n$), and the {\it piecewise linearity w.r.t the tuning parameter $\lambda$~\cite{Rosset-PeacewiseLin}} (to reduce computational effort by using a modification of the LARS algorithm).\\
\noindent We use the same notation as in Section~\ref{LassoOriginal} for the LASSO estimator. Set $ \hat{\beta}$ to be an
estimator of the regression vector $\beta$ based on $\mathbf{x}$ and $\mathbf{y}$. Let also $s$ be the sign vector of the
estimator $\hat{\beta}$. On the other hand, using the notation in Section~\ref{section-C-Lasso}, we set $\hat{\mu}=\widetilde{\mathbf{x}} \hat{\beta}$ where this time $\hat{\beta}$ is based on the augmented dataset $\widetilde{\mathbf{x}}$ and $\widetilde{\mathbf{y}}$.
\begin{ass}\label{betaGeneral}
The estimator $\hat{\mu}$ can be written as:
\begin{equation}
\hat{\mu} = u(\widetilde{\mathbf{x}},s) \widetilde{\mathbf{y}} + v(\widetilde{\mathbf{x}},s),
\end{equation}
where $u(\cdot)$ and $v(\cdot)$ are piecewise constant functions w.r.t.
$\widetilde{\mathbf{y}}$.
\end{ass}
As soon as Assumption~\ref{betaGeneral} holds, we can construct a conformal predictor corresponding to the estimator $\hat{\mu}$. Then many estimators can be considered. The CoLP and CoRP obviously belong to this class of predictors and we introduce here the Conformal Elastic Net Predictor (CENeP) which is a conformal predictor constructed based on the Elastic-Net modification of the LARS instead of the LASSO one ({\bf Step1} in Algorithm~\ref{alg:ConformalLasso}). This predictor
is defined by $u(\widetilde{\mathbf{x}},s) = \widetilde{\mathbf{x}}_{k}(\widetilde{\mathbf{x}}_{k}'\widetilde{\mathbf{x}}_{k}
+ \mu_{k}\mathbf{I}_{k} )^{-1}  \widetilde{\mathbf{x}}_{k}'$ and $v(\widetilde{\mathbf{x}},s) = -\lambda_{k} \widetilde{\mathbf{x}}_{k}(\widetilde{\mathbf{x}}_{k}'\widetilde{\mathbf{x}}_{k})^{-1} s_{k}$ where $\lambda_{k}$ and $\mu_{k}$ correspond respectively to the LASSO and Ridge tuning parameters in the definition of the Elastic-Net estimator and $\mathbf{I}_{k}$ is the $|\mathcal{A}_{k}| \times |\mathcal{A}_{k}|$ identity matrix~\cite{Zou-E-Net}. In the same way, we can define the Conformal Smooth Lasso Predictor (CoSmoLaP) based on a Smooth-Lasso modification of the LARS algorithm \cite{MoiSLasso}. Here $u(\widetilde{\mathbf{x}},s) = \widetilde{\mathbf{x}}_{k}(\widetilde{\mathbf{x}}_{k}'\widetilde{\mathbf{x}}_{k} + \mu_{k}\mathbf{J}_{k} )^{-1}  \widetilde{\mathbf{x}}_{k}'$ and $v(\widetilde{\mathbf{x}},s) = -\lambda_{k} \widetilde{\mathbf{x}}_{k}(\widetilde{\mathbf{x}}_{k}'\widetilde{\mathbf{x}}_{k})^{-1} s_{k}$. The difference between the CoSmoLaP definition the CENeP one is the identity matrix $\mathbf{I}_{k}$ which is replaced by the $|\mathcal{A}_{k}| \times
|\mathcal{A}_{k}|$ matrix $\mathbf{J}_{k}$ whose components are such that $(\mathbf{J}_{k})_{i,i} = 1$ if $i=1$ or $i= |\mathcal{A}_{k}|$ and $(\mathbf{J}_{k})_{i,i} = 2$ otherwise. Moreover for $(i,j)\in\{1,\ldots,\mathcal{A}_{k}\}^2$ with $i \neq j$, we have $(\mathbf{J}_{k})_{i,j} = -1$ if $|i-j| = 1$ and zero otherwise. Note that the definition of $\mathbf{J}_{k}$ makes the CoSmoLaP more appropriated to model with successive correlation between successive variables.

As for CoLP, we can define the nonconformity score of an expected label $y$ associated to the
estimator $\hat{\mu}$ as follows:
\begin{eqnarray*}
\left(\alpha_{1}(y),\ldots,\alpha_{n}(y)\right)' &:= &|\widetilde{\mathbf{y}}-\hat{\mu}| \\ &
= & |\left(\mathbf{I}-u(\widetilde{\mathbf{x}},s)\right)\widetilde{\mathbf{y}} - v(\widetilde{\mathbf{x}},s)|
\\ & = & |A + B\,y| ,
\end{eqnarray*}
with
\begin{equation*}
\left\{
\begin{array}{l}
A= (a_{1},\ldots,a_{n})' :=\left(\mathbf{I}-u(\widetilde{\mathbf{x}},s)\right) (y_{1},\ldots,y_{n-1},0)' - v(\widetilde{\mathbf{x}},s), \\
B= (b_{1},\ldots,b_{n})' :=\left(\mathbf{I}-u(\widetilde{\mathbf{x}},s)\right) (0,\ldots,0,1)', 
\end{array}
\right.
\end{equation*}
and $\mathbf{I}$ is the $n\times n$ identity matrix. The quantities $A$ and $B$ are the analogues of $A_{k}$ and $B_{k}$ respectively, when we considered the CoLP at the transition point $\lambda_{k},\, k = 1,\ldots,K$. Then replacing $A_{k}$ and $B_{k}$ by respectively $A$ and $B$ in {\bf Step~2.a} of Algorithm~\ref{alg:ConformalLasso}, we obtain the conformal predictors associated to the estimator $\hat{\mu}$.

Note that the dependency in the tuning parameter, noted $\lambda$, can be included in $u(\widetilde{\mathbf{x}},s)$ (as for CoRP) or $v(\widetilde{\mathbf{x}},s)$ or in both of them (as for the CoLP). For instance, in the construction of the CoLP,
this dependency is underlined in the matrix $\widetilde{\mathbf{x}}_k$ and the sign vector $s_k$ as they were computed by the LARS algorithm for a specified value $\lambda_{k}$ of the tuning parameter $\lambda$.

Computational cost of the construction of conformal predictors has also to be considered. Three main points interfere. First, one run of the LARS algorithm requires the same cost as the computation of the least square estimation. Then we have to consider the number of conformal predictors we have to construct: each value of the tuning parameter $\lambda$ provides a conformal predictor $\Gamma_{\lambda}$ using the algorithm described in Section~\ref{sec:Algo}. The final conformal predictor $\Gamma_{opt}$ is then the one with the minimal length. As for the CoRP, the main problem is: how many $\lambda$'s do we have to test? One way is to use a grid of value for $\lambda$ which lets open the problem of the choice of the grid and the window of this grid.\\
On the other hand, we saw how the LARS algorithm permits to reduce considerably the number of tuning parameters to be considered. Indeed the grid of tuning parameters values is directly described by the transition points $\lambda_{1},\ldots,\lambda_{K}$ obtained from the run of the LARS algorithm. Finally, let us consider {\it the construction of the conformal predictor itself:} this point has been treated in Vovk et al. \cite[Chapter 2.3 and 4.1]{Vovk_livre}. It turns out that sparse conformal predictors and the CoLP requires computation time $\mathcal{O}(n^{2})$ and can be reduced to $\mathcal{O}(n \log (n))$.
\section{Experimental Results}\label{sec:simul}
In the section we present the experimental performances of the Sparse Conformal Predictors (SCP) w.r.t. their validity, their accuracy and also their selection power. As benchmark, we use the CoRP\footnote{We construct the CoRP associated to same tuning parameters as the CoLP (i.e., the transition points $\lambda_{k}$ observed in Section~\ref{sec:Algo}). Note that the performance would not be inflected as conformal predictors according to this method are almost embedded and changes sensitively while the tuning parameter varies. See~\cite[page 39]{Vovk_livre} for more details.} for its validity and accuracy and the original LASSO and Elastic-Net estimators for their selection\footnote{We use a $\mathop{\rm BIC}$-type criterion to select the optimal tuning parameter. Such a criterion is adapted to variable selection.} power.\\
We consider three SCPs: the Conformal Lasso Predictor (CoLP was introduced in Sections~\ref{section-C-Lasso} and~\ref{sec:Algo}) and the Conformal Elastic Net Predictor (CENeP was described in Section~\ref{sec:generalization}). The last SCP called Conformal Ridge Lasso Predictor (CoRLaP) is a mix of the CoRP and the CoLP. To construct the CoRLaP, we use the variables selected by the LASSO modification of the LARS algorithm (\textbf{Step 1} in Algorithm~\ref{alg:ConformalLasso} described in Section~\ref{sec:Algo}). Then we use these variables to construct a CoRP. This conformal predictor can be seen as a restricted CoRP. All conformal predictors are constructed with confidence level $1-\varepsilon = 90\%$.
\subsection{Simulated Experiments}\label{subsec:Sim}
We consider four simulations from the linear regression model
\begin{equation*}
y=\mathbf{X}'\beta + \sigma\xi, \quad \xi\sim\mathcal{N}(0,1),\,
\mathbf{X}=(\mathbf{X}_1,\ldots,\mathbf{X}_{50})'\in \mathbb{R}^{50},
\end{equation*}
with $\beta\in \mathbb{R}^{50}$. Hence $p=50$ through the simulations. Noise level $\sigma$ and the sample size $n$ are let free. They will be specified during experiments.
\begin{description}
\item[{\it Example}~(a)]$[n/\sigma]$: {\it Very Sparse and Correlated.} Here
only $\beta_1$ is nonzero and equals $5$. Moreover, the design correlations matrix
$\Sigma$ is described by $\Sigma_{j,k} = \exp(-|j-k|)$
for $(j,k)\in\{15,\ldots,35\}^{2}$ and $\Sigma_{j,k} =\mathbb{I}(j=k)$
otherwise where $\mathbb{I}(\cdot)$ is the indicator function.\\
\item[{\it Example}~(b)]$[n/\sigma]$: {\it Sparse and Correlated.} The
correlations are defined as in Example~(a) and the regression vector is given by $\beta_j = -5+ 0.2j$ for $j=1,\ldots,5$; $\beta_j = 4+ 0.2j$ for $j=10,\ldots,25$ 
and zero otherwise.\\
\item[{\it Example}~(c)]$[n/\sigma]$: {\it Sparse and Highly correlated.} We have $\beta_{j} =5$ for $j\in\{1,\ldots,15\}$ and zero otherwise. We construct three groups of correlated variables: $\Sigma_{j,k}= 1$ when $(j,k)$ belongs to $\{1,\ldots,5\}^2$, $\{6,\ldots,10\}^2$ and $\{11,\ldots,15\}^2$; $\Sigma_{j,k}= 1$ for $(j,k)\in\{16,\ldots,p\}^2$ if $j=k$ and zero otherwise.\\
\item[{\it Example}~(d)]$[n/\sigma]$: {\it Non Sparse and correlated.} Here $\beta_{j} = 3+0.2j$ for $j\in\{1,\ldots,p\}$ and the correlations are described by $\Sigma_{j,k} = \exp(-|j-k|)$ for $(j,k)\in\{1,\ldots,p\}^{2}$.\\
\end{description}
We consider separately the three points of interest:  accuracy, validity and selection.\\

\begin{figure}[t]
\vskip -0.2in
\begin{center}
\includegraphics[height=2.3in] {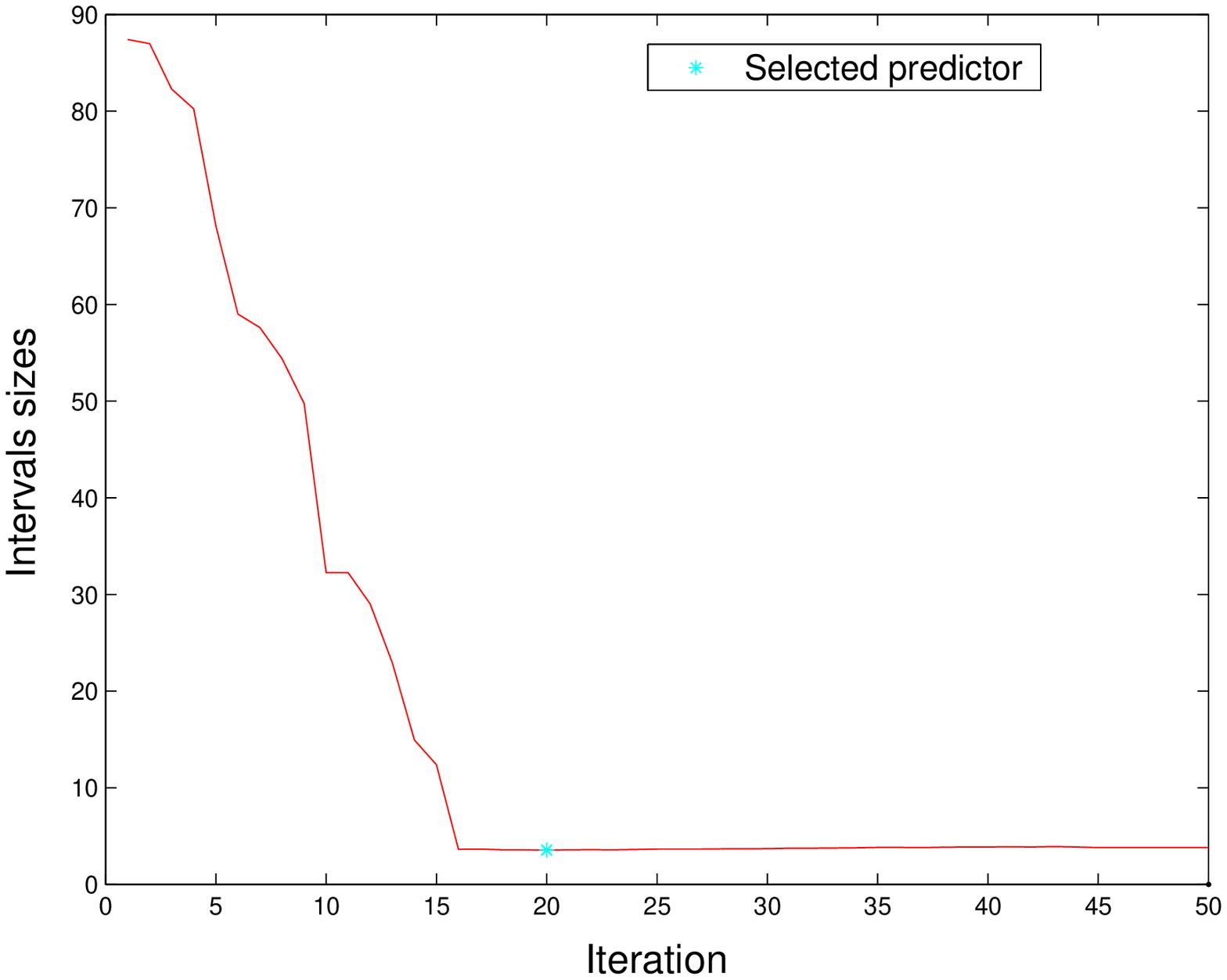}
\includegraphics[height=2.3in] {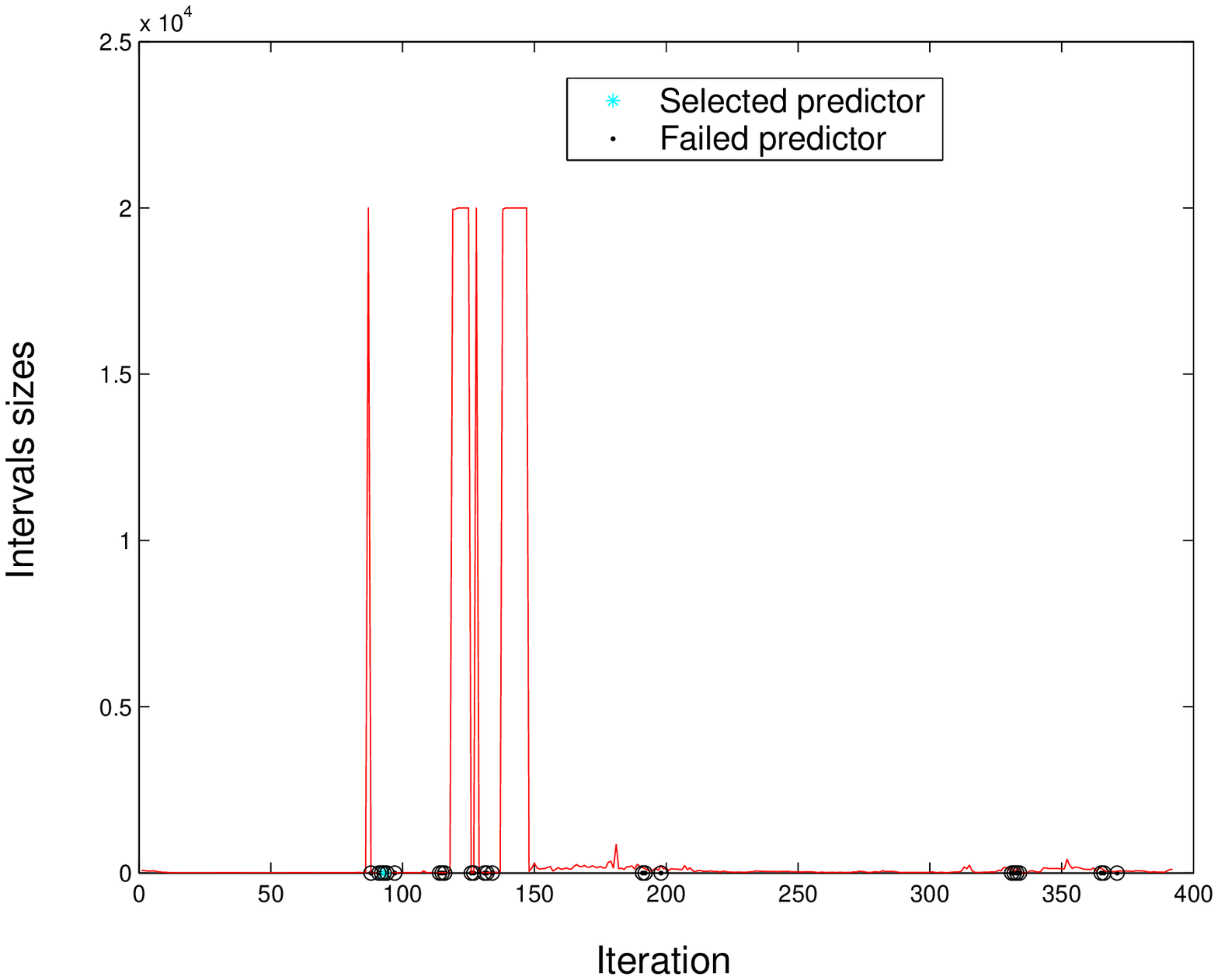}\\
\includegraphics[height=2in] {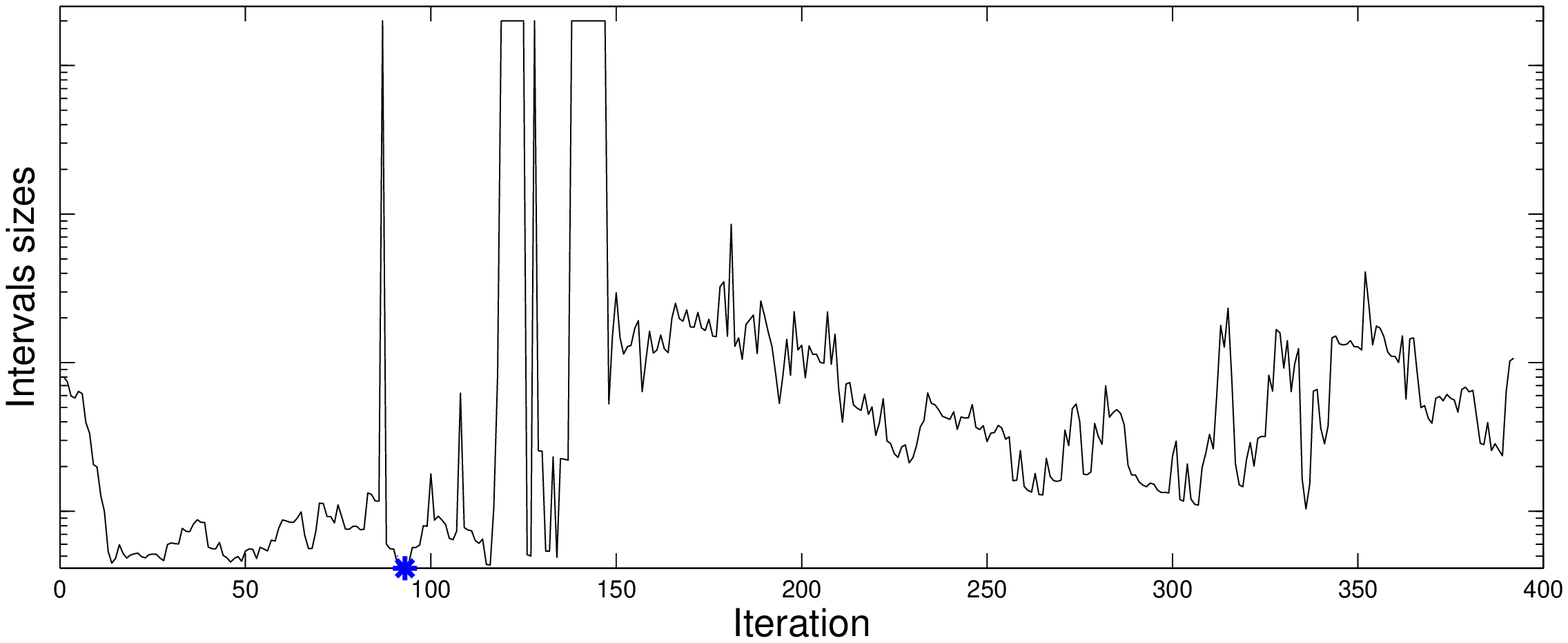}
\end{center}
\hfill\begin{minipage}[t]{0.90\textwidth}
\caption{\footnotesize Analysis of conformal predictors length (y-axis)
through the LASSO modification of the LARS algorithm iterations (x-axis: the first iteration corresponds to $\lambda_{max}$ and the last one corresponds to $\lambda_{min}$) in
Example~(c)[$300$/$1$] (top left) and in Example~(c)[$50$/$1$] (top right). The
iteration associated to the CoLP is marked by a blue star. Predictors which are
non valid are marked by a black circle. The panel of bottom shows the lengths of intervals in
a logarithmic scale.}
\end{minipage}
\label{fig:plotslength}
\end{figure}
\begin{description}
\item[Accuracy.] First of all, let us consider the length of the predictors $\Gamma_{k}^{\varepsilon},\, k=1,\ldots,K$ obtained at the end of \textbf{Step~2} in Algorithm~\ref{alg:ConformalLasso} described in Section~\ref{sec:Algo}. We remind that each of these predictors is associated to an iteration of a modification of the LARS algorithm, that is the transition points $\lambda_k ,\, k=1,\ldots,K$. Figure~\ref{fig:plotslength} illustrates the predictors lengths for the construction of the CoLP, when applied to Example~(c)[$n$/$1$] with $n=300$ and $n=50$. When $n=300$, we note that the length of the $\Gamma_{k}^{\varepsilon}$s sensitively changes from one iteration to the following and that the larger predictor has a reasonable length compared to the smallest one (about $10$ times larger). Then the construction is stable. We also observe that in the neighborhood of the optimal iteration (that is iteration $20$), the conformal predictors have approximately the same size. Such an observation can also be made when we take a look at Figure~\ref{fig:SuccesVsSteps} (left) when applied to Example~(b)[$300$/$1$]. On the other hand, when $n=50$, it appears that the predictors length grows drastically at some iteration (around iteration $85$). We even can not compare the lengths of the bigger and smaller predictors (more than $10^4$ times larger). In the same time, it seems that the construction becomes unstable as violent variations often happen after this iteration $85$. We will consider in the next point the validity of these predictors. However let us mention that in Example~(c)[$50$/$1$], the CoLP which is the smallest $\Gamma_{k}^{\varepsilon}$ and then the selected predictor is not valid (in Figure~\ref{fig:plotslength} (right), the selected predictor at iteration $93$ is not valid). This aspect can also be observed in Figure~\ref{fig:SuccesVsSteps} (right) (the graph corresponds to Example~(b)[$50$/$1$]) where the selected CoLP at iteration $57$ is not valid. Similar violent variations of the corresponding predictors lengths would have been observed after iteration $49$ if we have provided a graph as Figure~\ref{fig:plotslength} (right).\\
\begin{table}[t]
\caption{Validity frequencies [with precision $\pm 95\%$] of the CoRP, CoLP, CoRLaP, CENeP, the Early-Stopped CoLP and the $2$-PN~CoLP based on $1000$ replications.}
\label{tab:frqVSsigma} \vskip 0.15in
\begin{center}
\begin{sc}
\begin{tabular}{lc|cccr}
\hline
Example &  $\sigma$& CoRP & CoLP & CoRLaP & CENeP \\
\hline
(a)[$300$/$\sigma$] &  $1$    & 0.897$\pm$ 0.019& 0.876$\pm$ 0.020& 0.854$\pm$ 0.022& 0.878$\pm$ 0.020\\
&   $7$    & 0.894$\pm$ 0.019& 0.908$\pm$ 0.018& 0.894$\pm$ 0.019& 0.899$\pm$ 0.019 \\
&   $15$    & 0.893$\pm$ 0.019& 0.893 $\pm$ 0.019& 0.879$\pm$ 0.020& 0.887$\pm$ 0.020 \\
\hline
\hline
(b)[$300$/$\sigma$]&   $1$  & 0.901$\pm$ 0.018& 0.875$\pm$ 0.020& 0.869$\pm$ 0.021& 0.874$\pm$ 0.021\\
(c)[$300$/$\sigma$]&   $1$   & 0.900$\pm$ 0.019& 0.900$\pm$ 0.019& 0.891$\pm$ 0.019& 0.901$\pm$ 0.018 \\
(d)[$300$/$\sigma$]&   $1$    & 0.892$\pm$ 0.019& 0.895$\pm$ 0.019& 0.895$\pm$ 0.019& 0.895$\pm$ 0.019 \\
\hline
\hline
(a)[$50$/$\sigma$]&   $3$  & 0.887$\pm$ 0.020& 0.668$\pm$ 0.029& 0.414$\pm$ 0.030& 0.789$\pm$ 0.025\\
(a)[$20$/$\sigma$]&   $3$    & 0.865$\pm$ 0.021& 0.596$\pm$ 0.030& 0.304$\pm$ 0.028& 0.685$\pm$ 0.029 \\
\hline
\hline
Example &  $\sigma$& CoRP & CoLP & Stopped-CoLP & $2$-PN-CoLP \\
\hline
(a)[$50$/$\sigma$] &  $7$    & 0.853$\pm$ 0.022& 0.620$\pm$ 0.030 & 0.815$\pm$ 0.024& 0.881$\pm$ 0.020 \\
(b)[$50$/$\sigma$] &  $1$  & 0.875$\pm$ 0.020& 0.558$\pm$ 0.031& 0.814$\pm$ 0.024& 0.907 $\pm$ 0.018 \\
(c)[$20$/$\sigma$]&   $15$   & 0.875$\pm$ 0.020& 0.608$\pm$ 0.030& 0.769$\pm$ 0.026& 0.893$\pm$ 0.019 \\
(d)[$20$/$\sigma$]&   $1$   & 0.900$\pm$ 0.019 & 0.602$\pm$ 0.030& 0.793$\pm$ 0.025& 0.892$\pm$ 0.019 \\
\hline
\end{tabular}
\end{sc}
\end{center}
\vskip -0.1in
\end{table}
\item[Validity.] Now, we consider the validity of the selected predictors (cf. \textbf{Step 3} in Algorithm~\ref{alg:ConformalLasso}). As shown in Table~\ref{tab:frqVSsigma}, we observe that variations on the noise level, the variables correlations and the sparsity of the model do to not perturb the validity whereas the sample size relatively to the dimension $p$ does. When $n=300>p$, all the procedures seem to be quite similar and produce good predictors. In the other cases,
i.e., when $n=p=50$ and $n=20<p$, the selected confidence predictors have worst performance than expected (validity with smaller proportion than $1-\varepsilon = 90\%$). Moreover, Sparse Confidence Predictors perform worst than the CoRP as observed in Table~\ref{tab:frqVSsigma}. As pointed in the accuracy part, one explication can be observed in Figure~\ref{fig:plotslength} as the selected predictor which also is not valid (iteration $93$) corresponds to an iteration in the unstable zone (that is, after iteration $85$). Then in order to reduce the gap between SCP and CoRP in the cases $p\geq n$, we suggest to modify the selection criterion in \textbf{Step 3} in two ways. i) {\it Early Stopping CoLP:} do not consider (and do not construct) all the conformal predictors $\Gamma_{k}^{\varepsilon}$. Stop the construction of the predictors $\Gamma_{k}^{\varepsilon}$ as soon as the length of $\Gamma_{k}^{\varepsilon}$ (predictor at iteration $k$) has a length at least $10$ times larger than $\Gamma_{k-1}^{\varepsilon}$;
ii) {\it N Previous Neighbors CoLP:} we can enforce the Early Stopping rule by considering as final predictor: $\Gamma_{opt}^{\varepsilon} = \bigcup_{j:\,0\leq k-j<N} \Gamma_{j}^{\varepsilon}$, where $k$ is the index of the (selected) smallest predictor and $N$ is the number of neighbors we consider. Note that this method does not alter selection properties as $\Gamma_{k}^{\varepsilon}$ is usually constructed with more variables than $\Gamma_{j}^{\varepsilon},\, j < k$. It further does not alter a lot the accuracy as the Early Stopping rule ensures that we are in stable zone (cf. Figure~\ref{fig:plotslength} (right) and Figure~\ref{fig:SuccesVsSteps} (right)). Table~\ref{tab:frqVSsigma} sums up the performances of the early-stopped CoLP and the 2-PN CoLP in term of validity. We observe the good adaptation of both methods to the case $p = n$ and we remark that 2-PN CoLP nicely produce valid predictor even in the case $p>n$. This improvement in the term of validity can also be illustrated by Figure~\ref{fig:SuccesVsSteps} (right) where we observe that in Example~(b)[50/1], the early-stopped CoLP is valid whereas the original CoLP is not.\\
\begin{figure}[t]
\vskip -0.2in
\begin{center}
\includegraphics[height=2.3in] {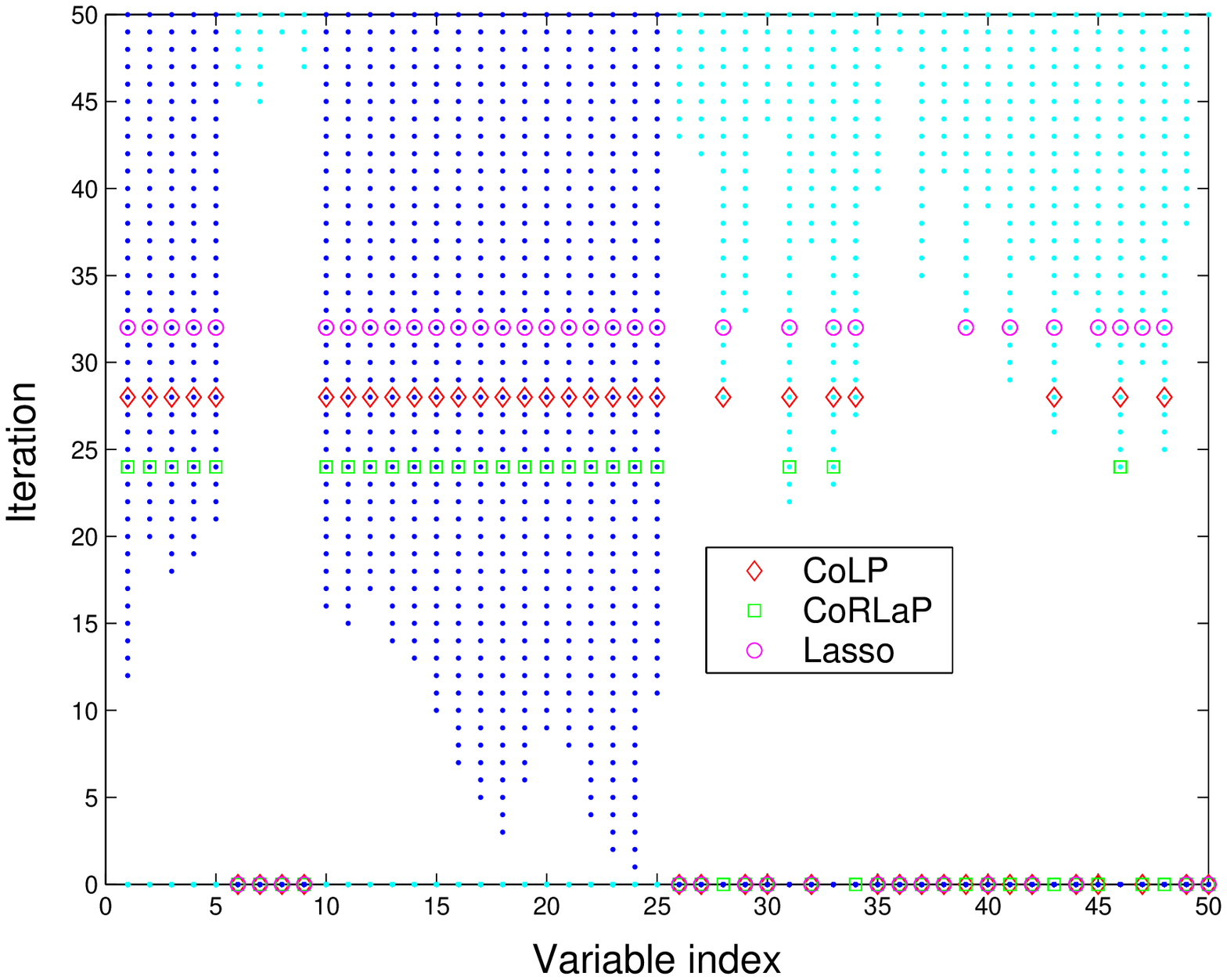}
\includegraphics[height=2.3in] {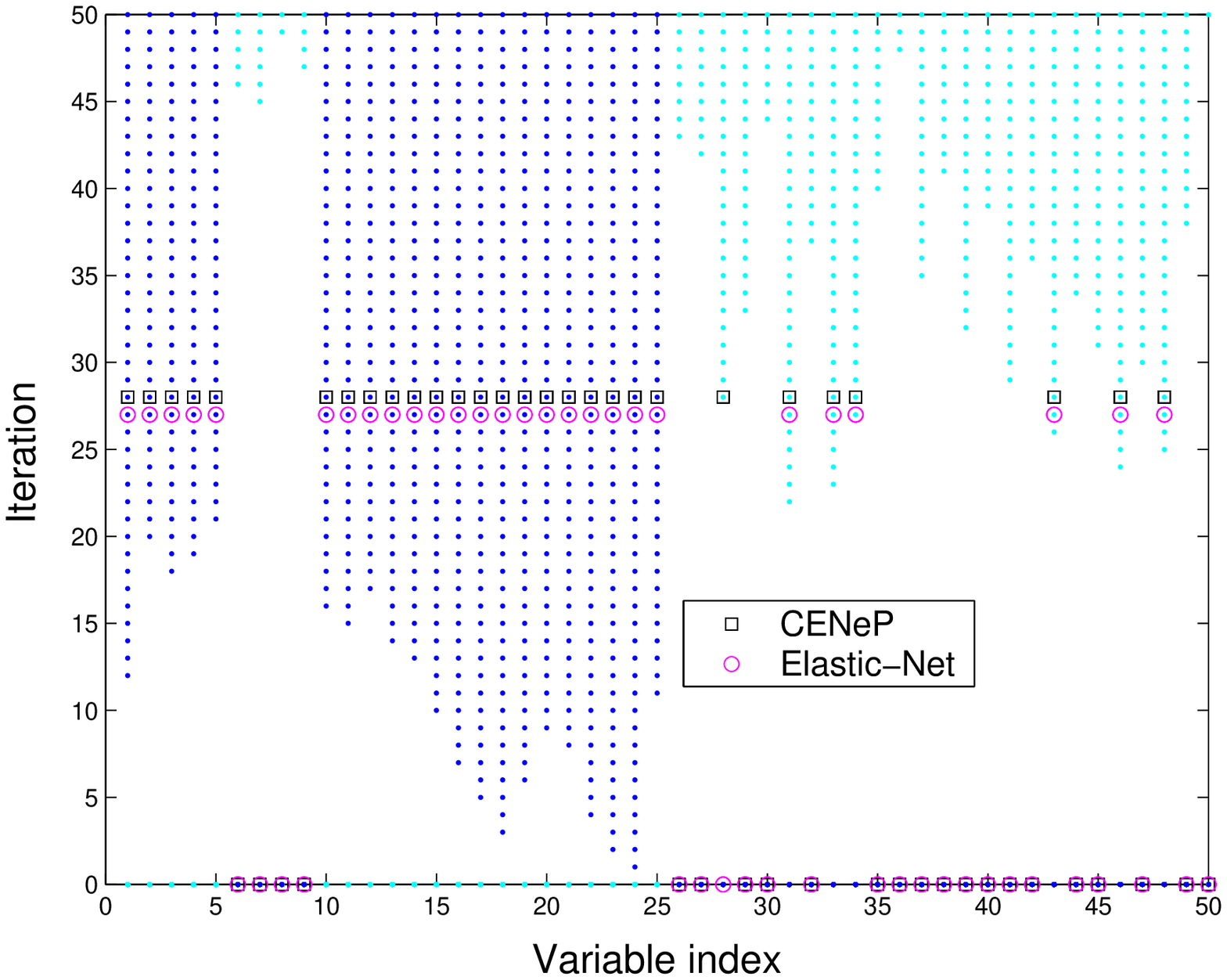}
\end{center}
\hfill\begin{minipage}[t]{0.90\textwidth}
\caption{\footnotesize Variable selection analysis for the CoLP, the CoRLaP and the CENeP in Example~(b)[300/1] (variables $1$ to $5$ and $10$ to $25$ are relevant; see variables in dark blue on the plot). On the left, we consider the CoLP and the CoRLaP selected variables (x-axis) with respect to the LASSO modification of the LARS algorithm iterations (y-axis: the first iteration corresponds to $\lambda_{max}$ and the last one corresponds to $\lambda_{min}$). On the right, we consider the CENeP selected variables (x-axis) with respect to the Elastic-Net modification of the LARS algorithm iterations (y-axis: the first iteration corresponds to $\lambda_{max}$ and the last one corresponds to $\lambda_{min}$). The selected iteration is marked by red diamonds for the CoLP, green squares for CoRLaP and black squares for the CENeP.}
\end{minipage}
\label{fig:plotsIndice}
\end{figure}
\item[Selection.] The selection ability of Sparse Conformal Predictors is here in concern. First, note that the selected variables in SCPs are directly linked to the selection ordering through the iterations of the LASSO or Elastic-Net modification of the LARS algorithm. Then, if the used modification of the LARS algorithm fails to recover the true model, we can not hope to get a predictor which contains only the true variables. Figure~\ref{fig:plotsIndice} illustrates the evolution of the variable selection of CoLP, CoRLaP and the LASSO on one hand and the CENeP and the Elastic-Net on the other hand, in Example~(b)[300/1]. It turns out that CoLP and CENeP select larger model that expected (that is, some noise variables are selected), as the LASSO and the Elastic-Net do. Moreover CoRLaP uses to select
a smaller subset of variables than the CoLP. Then it often produces a better variable selection performance than the other methods. It often provides closer model to the true one. Compared to the LASSO, it seems that the CoLP and the CoRLaP perform better in this example. However, we can not conclude the superiority of the CoLP on the LASSO in term of variable selection. A similar conclusion can be given when we compare the CENeP and the Elastic-Net. Nevertheless, the CENeP seems to select little larger models than the Elastic-Net. Finally, analogously to the superiority of the Elastic-Net compared to the LASSO, we can remark that the CENeP manages to have better selection performances compared to the CoLP and the CoRLaP when a group structure may exist between different variables (for instance in Example~(d)[$n$/$\sigma$]). This is due to the LASSO modification of the LARS algorithm which uses to select some noise variables before relevant ones in such cases.
\end{description}
\subsection{Real data}
We applied SCPs on $150$ randomly permutations of the House Boston dataset\footnote{The data and their description are available at http://archive.ics.uci.edu/ml/datasets/Housing.}, in which we randomly choose one row to be the new pair $(x_{new},y_{new})$. The original dataset consists of $506$ observations with $13$ variables. When we consider variable selection, we note that almost all SCPs are constructed without the variable $X_7=(x_{1,7},\ldots,x_{505,7})$. This variable is selected with frequencies lower than $3\%$. The CoRLaP also does not consider the variable $X_3$ as relevant with a frequence equal to $17\%$. Conforming to Section~\ref{subsec:Sim}, we would better consider $X_3$ irrelevant as the CoRLaP uses to produce better performance when variable selection is in concern. Then we conclude that the proportion of non-retail business acres per town and the proportion of owner-occupied units built prior to 1940 do not interfere in the value of owner-occupied homes. We also can notice that variable selection sligtly improved accuracy of conformal predictors in all presented experiments. Here, we can for instance remark that the median lengths of the CoLP, the CoRLaP and the CENeP are respectively $13.61$, $13.50$ and $13.58$, whereas CoRP length is $14.45$.

\section{Conclusion}
We presented Sparse Conformal Predictors, a family of $l_1$ regularized conformal predictors. We focused on LASSO and Elastic-Net versions of these Sparse Conformal Predictors. We illustrated their performance in term of accuracy, validity and variable selection. We concluded that such Sparse Conformal Predictors are valid and nicely exploit the sparsity of the model when the sample size is larger than the the number of variables (i.e, when $n>p$). We also provided a way to adopt these sparse predictors to the case $p\geq n$ through a pair of rules we called Early Stopping and $N$ Previous Neighbors rules.\\
Several extensions of this work can be explored such as the construction of SCP with Adaptive LASSO~\cite{Zou-Adaptive-Lasso} and they will be investigated in future work.
\begin{flushright}
$ \Box $
 \end{flushright}
\noindent {\bf Acknowledgement.} We would like to thank Professor Arnak Dalalyan and Professor Nicolas Vayatis for insightful comments.

\bibliographystyle{plain}

\bibliography{Conformal_lasso}

\end{document}